\documentclass{gtart}
\def\ifplaintex{\expandafter\ifx\csname documentclass\endcsname\relax}

\def\gtp{{\mathsurround=0pt\it $\cal G\mskip-2mu$eometry \&\ 
$\cal T\!\!$opology $\cal P\!$ublications}}  % GT publications

\def\recd{{\small Received:\qua\receiveddate\ifx\reviseddate\relax
\else\qquad Revised:\qua\reviseddate\fi\par}} 

\def\lognumber#1{\def\thelognumber{#1}}
\def\volumenumber#1{\def\thevolumenumber{#1}}
\def\volumeyear#1{\def\thevolumeyear{#1}}
\def\papernumber#1{\def\thepapernumber{#1}}
\def\pagenumbers#1#2{\def\startpage{#1}\def\finishpage{#2}}
\def\published#1{\def\publishdate{#1}}

\def\received#1{\def\receiveddate{#1}}
\def\revised#1{\def\reviseddate{#1}}
\def\accepted#1{\def\accepteddate{#1}}

\long\def\asciiabstract#1{\long\def\theasciiabstract{#1}}

\let\\\par\let\thelognumber\relax\let\thevolumenumber\relax
\let\thepapernumber\relax\let\thevolumeyear\relax\let\startpage\relax
\let\finishpage\relax\let\publishdate\relax\let\receiveddate\relax
\let\reviseddate\relax\let\accepteddate\relax\let\theasciititle\relax
\let\theasciiauthors\relax
\let\theasciiabstract\relax

\let\theasciiemail\relax

\ifplaintex
\font\logobig=cmssbx10 scaled 3836
\font\logomed=cmssbx10 scaled 2557
\else
\font\logobig=cmssbx10 scaled 4200
\font\logomed=cmssbx10 scaled 2800
\fi

\long\def\makeagttitle{   %%% start of definition of \makeagttitle
\count0=\startpage
\agt\hfill      %   Journal title (top left) 
\hbox to 45truept{\vbox to 0pt{\vglue -13truept{\logomed A\kern -.37em{\logobig 
T}\kern -.38em G}\vss}\hss}
\break
{\small Volume \thevolumenumber\ (\thevolumeyear)
\startpage--\finishpage\nl
Published: \publishdate}

\vglue .25truein

{\parskip=0pt\leftskip 0pt plus
1fil\def\\{\par\smallskip}{\Large\bf\thetitle}\par\medskip} \vglue
0.05truein

{\parskip=0pt\leftskip 0pt plus 1fil\def\\{\par}{\sc\theauthors}
\par\medskip}%
 
\vglue 0.03truein 

{\small\leftskip 25truept\rightskip 25truept{\bf Abstract}\stdspace\theabstract

{\bf AMS Classification}\stdspace\theprimaryclass
\ifx\thesecondaryclass\relax\else; \thesecondaryclass\fi\par
{\bf Keywords}\stdspace \thekeywords\par}\vglue 7truept

}   %%%% end of definition of \makeagttitle

\ifplaintex
\font\phead=cmsl9 scaled 950
\font\pnum=cmbx10 scaled 913
\font\pfoot=cmsl9 scaled 950
\headline{\vbox to 0pt{\vskip -4.5mm\line{\small\phead\ifnum
\count0=\startpage ISSN 1472-2739 (on-line) 1472-2747 (printed)
\hfill {\pnum\folio}\else\ifodd\count0\def\\{ }% 
\ifx\theshorttitle\relax\thetitle\else\theshorttitle\fi\hfill{\pnum\folio}
\else\def\\{ and }{\pnum\folio}\hfill\ifx\theshortauthors\relax\theauthors
\else\theshortauthors\fi\fi\fi}\vss}}
\footline{\vbox to 0pt{\vglue 0mm\line{\small\pfoot\ifnum\count0=\startpage
\copyright\ \gtp\hfill\else
\agt, Volume \thevolumenumber\ (\thevolumeyear)\hfill\fi}\vss}}
\else
\font=cmsl9 scaled 1050
\font\lnum=cmbx10 
\font=cmsl9 scaled 1050
\makeatletter
\def\@oddhead{{\small\lhead\ifnum\count0=\startpage ISSN 1472-2739 
(on-line) 1472-2747 (printed)\hfill {\lnum\number\count0}\else\ifodd\count0
\def\\{ }\ifx\theshorttitle\relax \thetitle \else\theshorttitle\fi\hfill
{\lnum\number\count0}\else\def\\{ and }{\lnum\number\count0}
\hfill\ifx\theshortauthors\relax 
\theauthors\else\theshortauthors\fi\fi\fi}}\def\@evenhead{\@oddhead}
\def\@oddfoot{\small\lfoot\ifnum\count0=\startpage\copyright\ \gtp\hfill\else
\agt, Volume \thevolumenumber\ (\thevolumeyear)\hfill\fi}
\def\@evenfoot{\@oddfoot}
\makeatother
\fi
\let\maketitlepage\makeagttitle
\let\makeshorttitle\maketitlepage
\let\maketitle\maketitlepage

\newwrite\gtoutfile
\long\gdef\makeheadfile{  %%% start of definition of \makeheadfile
{\def\\{, }\def\s{ }
\immediate\openout\gtoutfile head.xxx
\immediate\write\gtoutfile{To: math@arxiv.org}
\immediate\write\gtoutfile{Subject: put OR rep NNNNN:ppppp}
\immediate\write\gtoutfile{--text follows this line--}
\immediate\write\gtoutfile{Proxy-for: \ifx\theasciiauthors\relax
\theauthors\else\theasciiauthors\fi\s<\ifx\theasciiemail\relax\theemail\else\theasciiemail\fi>}
\immediate\write\gtoutfile{\noexpand\\}
\immediate\write\gtoutfile{Authors: \ifx\theasciiauthors\relax
\theauthors\else\theasciiauthors\fi}
{\def\\{ }\immediate\write\gtoutfile{Title: \ifx\theasciititle\relax
\thetitle\else\theasciititle\fi}}
\immediate\write\gtoutfile{Subj-class: GT or SG, GR etc}
\immediate\write\gtoutfile{MSC-class: \theprimaryclass\ifx\thesecondaryclass\relax\else, \thesecondaryclass\fi}
\immediate\write\gtoutfile{Journal-ref: Algebr. Geom. Topol. \thevolumenumber\s
(\thevolumeyear) \startpage-\finishpage}
\immediate\write\gtoutfile{Comments: Published by Algebraic and
Geometric Topology at}
\immediate\write\gtoutfile{\s\s\s  http://www.maths.warwick.ac.uk/agt/AGTVol\thevolumenumber/agt-\thevolumenumber-\thepapernumber.abs.html}
\immediate\write\gtoutfile{\noexpand\\}
\immediate\write\gtoutfile{}
\ifx\theasciiabstract\relax
\immediate\write\gtoutfile{\theabstract}\else
\immediate\write\gtoutfile{\theasciiabstract}\fi
\immediate\write\gtoutfile{}
\immediate\write\gtoutfile{\noexpand\\}
\immediate\write\gtoutfile{}
\immediate\closeout\gtoutfile}}  %%% end of definition of \makeheadfile

\def\maketitlepage{\makeagttitle\makeheadfile}
\let\makeshorttitle\maketitlepage
\let\maketitle\maketitlepage

\def\ifplaintex{\expandafter\ifx\csname documentclass\endcsname\relax}

\def\gtp{{\mathsurround=0pt\it $\cal G\mskip-2mu$eometry \&\ 
$\cal T\!\!$opology $\cal P\!$ublications}}  % GT publications

\def\recd{{\small Received:\qua\receiveddate\ifx\reviseddate\relax
\else\qquad Revised:\qua\reviseddate\fi\par}} 

\def\lognumber#1{\def\thelognumber{#1}}
\def\volumenumber#1{\def\thevolumenumber{#1}}
\def\volumeyear#1{\def\thevolumeyear{#1}}
\def\papernumber#1{\def\thepapernumber{#1}}
\def\pagenumbers#1#2{\def\startpage{#1}\def\finishpage{#2}}
\def\published#1{\def\publishdate{#1}}

\def\received#1{\def\receiveddate{#1}}
\def\revised#1{\def\reviseddate{#1}}
\def\accepted#1{\def\accepteddate{#1}}

\long\def\asciiabstract#1{\long\def\theasciiabstract{#1}}

\let\\\par\let\thelognumber\relax\let\thevolumenumber\relax
\let\thepapernumber\relax\let\thevolumeyear\relax\let\startpage\relax
\let\finishpage\relax\let\publishdate\relax\let\receiveddate\relax
\let\reviseddate\relax\let\accepteddate\relax\let\theasciititle\relax
\let\theasciiauthors\relax
\let\theasciiabstract\relax

\let\theasciiemail\relax

\ifplaintex
\font\logobig=cmssbx10 scaled 3836
\font\logomed=cmssbx10 scaled 2557
\else
\font\logobig=cmssbx10 scaled 4200
\font\logomed=cmssbx10 scaled 2800
\fi

\long\def\makeagttitle{   %%% start of definition of \makeagttitle
\count0=\startpage
\agt\hfill      %   Journal title (top left) 
\hbox to 45truept{\vbox to 0pt{\vglue -13truept{\logomed A\kern -.37em{\logobig 
T}\kern -.38em G}\vss}\hss}
\break
{\small Volume \thevolumenumber\ (\thevolumeyear)
\startpage--\finishpage\nl
Published: \publishdate}

\vglue .25truein

{\parskip=0pt\leftskip 0pt plus
1fil\def\\{\par\smallskip}{\Large\bf\thetitle}\par\medskip} \vglue
0.05truein

{\parskip=0pt\leftskip 0pt plus 1fil\def\\{\par}{\sc\theauthors}
\par\medskip}%
 
\vglue 0.03truein 

{\small\leftskip 25truept\rightskip 25truept{\bf Abstract}\stdspace\theabstract

{\bf AMS Classification}\stdspace\theprimaryclass
\ifx\thesecondaryclass\relax\else; \thesecondaryclass\fi\par
{\bf Keywords}\stdspace \thekeywords\par}\vglue 7truept

}   %%%% end of definition of \makeagttitle

\ifplaintex
\font\phead=cmsl9 scaled 950
\font\pnum=cmbx10 scaled 913
\font\pfoot=cmsl9 scaled 950
\headline{\vbox to 0pt{\vskip -4.5mm\line{\small\phead\ifnum
\count0=\startpage ISSN 1472-2739 (on-line) 1472-2747 (printed)
\hfill {\pnum\folio}\else\ifodd\count0\def\\{ }% 
\ifx\theshorttitle\relax\thetitle\else\theshorttitle\fi\hfill{\pnum\folio}
\else\def\\{ and }{\pnum\folio}\hfill\ifx\theshortauthors\relax\theauthors
\else\theshortauthors\fi\fi\fi}\vss}}
\footline{\vbox to 0pt{\vglue 0mm\line{\small\pfoot\ifnum\count0=\startpage
\copyright\ \gtp\hfill\else
\agt, Volume \thevolumenumber\ (\thevolumeyear)\hfill\fi}\vss}}
\else
\font=cmsl9 scaled 1050
\font\lnum=cmbx10 
\font=cmsl9 scaled 1050
\makeatletter
\def\@oddhead{{\small\lhead\ifnum\count0=\startpage ISSN 1472-2739 
(on-line) 1472-2747 (printed)\hfill {\lnum\number\count0}\else\ifodd\count0
\def\\{ }\ifx\theshorttitle\relax \thetitle \else\theshorttitle\fi\hfill
{\lnum\number\count0}\else\def\\{ and }{\lnum\number\count0}
\hfill\ifx\theshortauthors\relax 
\theauthors\else\theshortauthors\fi\fi\fi}}\def\@evenhead{\@oddhead}
\def\@oddfoot{\small\lfoot\ifnum\count0=\startpage\copyright\ \gtp\hfill\else
\agt, Volume \thevolumenumber\ (\thevolumeyear)\hfill\fi}
\def\@evenfoot{\@oddfoot}
\makeatother
\fi
\let\maketitlepage\makeagttitle
\let\makeshorttitle\maketitlepage
\let\maketitle\maketitlepage

\newwrite\gtoutfile
\long\gdef\makeheadfile{  %%% start of definition of \makeheadfile
{\def\\{, }\def\s{ }
\immediate\openout\gtoutfile head.xxx
\immediate\write\gtoutfile{To: math@arxiv.org}
\immediate\write\gtoutfile{Subject: put OR rep NNNNN:ppppp}
\immediate\write\gtoutfile{--text follows this line--}
\immediate\write\gtoutfile{Proxy-for: \ifx\theasciiauthors\relax
\theauthors\else\theasciiauthors\fi\s<\ifx\theasciiemail\relax\theemail\else\theasciiemail\fi>}
\immediate\write\gtoutfile{\noexpand\\}
\immediate\write\gtoutfile{Authors: \ifx\theasciiauthors\relax
\theauthors\else\theasciiauthors\fi}
{\def\\{ }\immediate\write\gtoutfile{Title: \ifx\theasciititle\relax
\thetitle\else\theasciititle\fi}}
\immediate\write\gtoutfile{Subj-class: GT or SG, GR etc}
\immediate\write\gtoutfile{MSC-class: \theprimaryclass\ifx\thesecondaryclass\relax\else, \thesecondaryclass\fi}
\immediate\write\gtoutfile{Journal-ref: Algebr. Geom. Topol. \thevolumenumber\s
(\thevolumeyear) \startpage-\finishpage}
\immediate\write\gtoutfile{Comments: Published by Algebraic and
Geometric Topology at}
\immediate\write\gtoutfile{\s\s\s  http://www.maths.warwick.ac.uk/agt/AGTVol\thevolumenumber/agt-\thevolumenumber-\thepapernumber.abs.html}
\immediate\write\gtoutfile{\noexpand\\}
\immediate\write\gtoutfile{}
\ifx\theasciiabstract\relax
\immediate\write\gtoutfile{\theabstract}\else
\immediate\write\gtoutfile{\theasciiabstract}\fi
\immediate\write\gtoutfile{}
\immediate\write\gtoutfile{\noexpand\\}
\immediate\write\gtoutfile{}
\immediate\closeout\gtoutfile}}  %%% end of definition of \makeheadfile

\def\maketitlepage{\makeagttitle\makeheadfile}
\let\makeshorttitle\maketitlepage
\let\maketitle\maketitlepage

\lognumber{21}
\volumenumber{1}
\volumeyear{2001}
\papernumber{21}
\pagenumbers{427}{434}
\received{24 May 2001}
\revised{26 July 2001}
\accepted{27 July 2001}
\published{31 July 2001}

\usepackage{graphicx,amsfonts,latexsym,amsmath}

\input prepictex
\input pictex
\input postpictex

\newtheorem{theorem}{Theorem}
\newtheorem{lemma}[theorem]{Lemma}

\newcommand{\R}{\mathbb{R}}
\newcommand{\Z}{\mathbb{Z}}
\newcommand{\tb}{\overline{tb}}
\newcommand{\mtb}[1]{\overline{tb}\:\!(#1)}
\newcommand{\mindeg}{\textrm{min-deg}_a}

\begin{document}

\title{Maximal Thurston-Bennequin Number\\of Two-Bridge Links}
\author{Lenhard L. Ng}
\address{Department of Mathematics, Massachusetts Institute of Technology,\\
77 Massachusetts Avenue, Cambridge, MA 02139, USA}
\email{ng@alum.mit.edu}
\url{http://www-math.mit.edu/\char'176lenny/}

\begin{abstract}
We compute the maximal Thurston-Bennequin number for a Legendrian
two-bridge knot or oriented two-bridge link in standard contact $\R^3$,
by showing that the upper bound given by the Kauffman polynomial is sharp.
As an application, we present a table of maximal Thurston-Bennequin numbers
for prime knots with nine or fewer crossings.
\end{abstract}

\asciiabstract{
We compute the maximal Thurston-Bennequin number for a Legendrian
two-bridge knot or oriented two-bridge link in standard contact R^3,
by showing that the upper bound given by the Kauffman polynomial is sharp.
As an application, we present a table of maximal Thurston-Bennequin numbers
for prime knots with nine or fewer crossings.}

\primaryclass{53D12}
\secondaryclass{57M15}
\keywords{Legendrian knot, two-bridge, 
Thurston-Bennequin number, Kauffman polynomial}

\maketitle

\section{Introduction}

A {\it Legendrian knot or link} in standard contact $\R^3$ is a knot or link 
which is everywhere tangent to the two-plane distribution induced by the 
contact one-form $dz - y\,dx$.  Given either a Legendrian knot or an 
oriented Legendrian link, we may define its
Thurston-Bennequin number, abbreviated $tb$, which is a Legendrian 
isotopy invariant; see, e.g., \cite{Ben} or \cite{FT}.  
(Henceforth, we will use the word ``link'' to denote
either a knot or an oriented link.)  For a fixed smooth link type $K$, 
the set of possible
Thurston-Bennequin numbers of Legendrian links in $\R^3$ isotopic to $K$
is bounded above; it is then natural to try to compute the maximum
$\mtb{K}$ of $tb$ over all such links.
Note that we distinguish between a link and its mirror; $\tb$ is
often different for the two.

Bennequin \cite{Ben} proved the first upper bound on $\mtb{K}$, 
in terms of the (three-ball) genus of $K$.  Since then, other
upper bounds have been found in terms of the HOMFLY and Kauffman polynomials
of $K$.  The strongest upper bound, in general, seems to be the Kauffman
bound first discovered by Rudolph \cite{Rud}, with alternative proofs
given by several authors;
see \cite{Fer} for a more detailed history of the subject.

Let $F_K(a,x)$ be the Kauffman polynomial of a link $K$,
and let $\mindeg$ denote the minimum degree in the framing variable $a$.
With the normalizations of \cite{FT}, the Kauffman bound states that
\[
\mtb{K} \leq \mindeg F_K(a,x) - 1.
\]

The Kauffman inequality is not sharp in general; see, e.g., 
\cite{EH,Fer}.
Sharpness has been established, however, for some small classes of knots,
including positive knots \cite{Tan}, most torus knots
\cite{Eps,EH}, and most three-stranded pretzel knots
\cite{Ng}.  In this note, we will establish sharpness for a somewhat
``larger'' class of links, the $2$-bridge (rational) links.  (We remark that
the HOMFLY bound is not sharp in general for this class.)

\begin{theorem}
If $K$ is a $2$-bridge link, then $\mtb{K} = \mindeg F_K(a,x) - 1$.
\label{main}
\end{theorem}

Theorem~\ref{main} will be proved in Section~\ref{sec:proof}.

Recall that a $2$-bridge link is any nontrivial link which admits
a diagram with four vertical tangencies (two on the left, two on the right).
This class of links includes many prime knots with a small number of 
crossings.  More precisely, all prime knots with seven or fewer crossings
are $2$-bridge, as are all prime knots with eight or nine crossings
except the following: $8_5$, $8_{10}$,
$8_{15}$--$8_{21}$, $9_{16}$, $9_{22}$, $9_{24}$, $9_{25}$, $9_{28}$,
$9_{29}$, $9_{30}$, and $9_{32}$--$9_{49}$.

Hand-drawn examples by N.\ Yufa and the author \cite{Yu} show that the
Kauffman bound is sharp for all of the above non-$2$-bridge $8$-crossing 
knots, except for $8_{19}$ (more precisely, 
the mirror image of the version drawn in \cite{Rol}).  Since $8_{19}$
is the $(4,-3)$ torus knot, a result of \cite{EH}
yields $\tb = -12$ in this case, while the Kauffman bound gives
$\tb \leq -11$.
Inspection of the non-prime knots with eight or fewer crossings
shows that the Kauffman bound is sharp for all such knots.
We thus have the following result.

\begin{theorem}
The Kauffman bound is sharp for all knots with eight or fewer
crossings, except the $(4,-3)$ torus knot $8_{19}$.
\label{small}
\end{theorem}

Further drawings show that the Kauffman bound is sharp for
all of the $9$-crossing prime knots which are not $2$-bridge, except
possibly for $9_{42}$ (more precisely, the mirror of the $9_{42}$ diagram in
\cite{Rol}).  For this last knot, we believe that $\tb = -5$, while
Kauffman gives $\tb \leq -3$.

An appendix to this note provides a table of
$\tb$ for prime knots with nine or fewer crossings.
Note that this table improves on the one from \cite{Tan}, 
which only considers one knot out of each mirror pair, and which does not
achieve sharpness in a number of cases.

{\bf Acknowledgements}
The author would like to thank Nataliya Yufa for her careful drawings which
culminated in Theorem~\ref{small} and the appendix table, Tom Mrowka 
and John Etnyre for useful discussions, and Isadore Singer for his 
encouragement and support.  This work was partially supported by
grants from the NSF and DOE.

\section{Proof of Theorem \ref{main}}
\label{sec:proof}

Let $K$ be a $2$-bridge link; we first need to find a suitable 
Legendrian embedding of $K$.  
Say that a link diagram is in {\it rational form} if it is in the form 
$T(a_1,\ldots,a_n)$ illustrated by Figure~\ref{bridge} for some
$a_1,\dots,a_n$.  Clearly any rational-form diagram corresponds to 
either the trivial knot or a $2$-bridge link;
by the classification of $2$-bridge links \cite{Sch}, 
any $2$-bridge link has a rational-form diagram.

\begin{figure}[ht!]\small
\centerline{
\font\thinlinefont=cmr5
\beginpicture
\setcoordinatesystem units <1.2cm,1.2cm>
\unitlength=1.00000cm
\linethickness=1pt
\setplotsymbol ({\makebox(0,0)[l]{\tencirc\symbol{'160}}})
\setshadesymbol ({\thinlinefont .})
\setlinear
\linethickness=1pt
\setplotsymbol ({\makebox(0,0)[l]{\tencirc\symbol{'160}}})
%
% Fig INTERPOLATED PT SPLINE
%
\plot  3.812 18.076 	 3.888 18.092
	 3.953 18.103
	 4.060 18.111
	 4.140 18.103
	 4.204 18.076
	 4.276 17.997
	 4.328 17.887
	 4.365 17.773
	 4.398 17.685
	 4.433 17.586
	 4.452 17.521
	 4.473 17.455
	 4.525 17.342
	 4.595 17.293
	 4.646 17.335
	 4.669 17.396
	 4.693 17.488
	/
\linethickness=1pt
\setplotsymbol ({\makebox(0,0)[l]{\tencirc\symbol{'160}}})
%
% Fig INTERPOLATED PT SPLINE
%
\plot  3.812 17.293 	 3.888 17.277
	 3.953 17.266
	 4.060 17.257
	 4.140 17.266
	 4.204 17.293
	 4.257 17.356
	 4.279 17.411
	 4.301 17.488
	/
\linethickness=1pt
\setplotsymbol ({\makebox(0,0)[l]{\tencirc\symbol{'160}}})
%
% Fig INTERPOLATED PT SPLINE
%
\plot  4.496 17.879 	 4.520 17.972
	 4.543 18.033
	 4.595 18.076
	 4.664 18.027
	 4.716 17.913
	 4.737 17.847
	 4.756 17.783
	 4.790 17.685
	 4.824 17.586
	 4.843 17.521
	 4.864 17.455
	 4.916 17.342
	 4.985 17.293
	 5.037 17.335
	 5.060 17.396
	 5.084 17.488
	/
\linethickness=1pt
\setplotsymbol ({\makebox(0,0)[l]{\tencirc\symbol{'160}}})
%
% Fig INTERPOLATED PT SPLINE
%
\plot  4.887 17.879 	 4.911 17.972
	 4.934 18.033
	 4.985 18.076
	 5.055 18.027
	 5.107 17.913
	 5.128 17.847
	 5.147 17.783
	 5.182 17.685
	 5.216 17.586
	 5.235 17.521
	 5.256 17.455
	 5.307 17.342
	 5.376 17.293
	 5.429 17.335
	 5.452 17.396
	 5.476 17.488
	/
\linethickness=1pt
\setplotsymbol ({\makebox(0,0)[l]{\tencirc\symbol{'160}}})
%
% Fig INTERPOLATED PT SPLINE
%
\plot  6.159 17.293 	 6.084 17.278
	 6.018 17.266
	 5.912 17.258
	 5.831 17.266
	 5.768 17.293
	 5.695 17.372
	 5.644 17.482
	 5.606 17.596
	 5.573 17.685
	 5.539 17.783
	 5.520 17.847
	 5.499 17.913
	 5.447 18.027
	 5.376 18.076
	 5.326 18.032
	 5.303 17.972
	 5.279 17.879
	/
\linethickness=1pt
\setplotsymbol ({\makebox(0,0)[l]{\tencirc\symbol{'160}}})
%
% Fig INTERPOLATED PT SPLINE
%
\plot  6.159 18.076 	 6.084 18.091
	 6.018 18.102
	 5.912 18.111
	 5.831 18.102
	 5.768 18.076
	 5.715 18.013
	 5.693 17.958
	 5.671 17.879
	/
\linethickness=1pt
\setplotsymbol ({\makebox(0,0)[l]{\tencirc\symbol{'160}}})
%
% Fig INTERPOLATED PT SPLINE
%
\plot  9.095 17.293 	 9.171 17.277
	 9.236 17.266
	 9.343 17.257
	 9.423 17.266
	 9.487 17.293
	 9.559 17.372
	 9.610 17.482
	 9.648 17.596
	 9.682 17.685
	 9.716 17.783
	 9.735 17.847
	 9.756 17.913
	 9.807 18.027
	 9.876 18.076
	 9.929 18.033
	 9.952 17.972
	 9.976 17.879
	/
\linethickness=1pt
\setplotsymbol ({\makebox(0,0)[l]{\tencirc\symbol{'160}}})
%
% Fig INTERPOLATED PT SPLINE
%
\plot  9.095 18.076 	 9.171 18.092
	 9.236 18.103
	 9.343 18.111
	 9.423 18.103
	 9.487 18.076
	 9.540 18.013
	 9.562 17.957
	 9.584 17.879
	/
\linethickness=1pt
\setplotsymbol ({\makebox(0,0)[l]{\tencirc\symbol{'160}}})
%
% Fig INTERPOLATED PT SPLINE
%
\plot  9.779 17.488 	 9.803 17.396
	 9.825 17.335
	 9.876 17.293
	 9.947 17.342
	 9.999 17.455
	10.020 17.521
	10.039 17.586
	10.073 17.685
	10.107 17.783
	10.126 17.847
	10.147 17.913
	10.199 18.027
	10.268 18.076
	10.320 18.033
	10.343 17.972
	10.367 17.879
	/
\linethickness=1pt
\setplotsymbol ({\makebox(0,0)[l]{\tencirc\symbol{'160}}})
%
% Fig INTERPOLATED PT SPLINE
%
\plot 10.171 17.488 	10.194 17.396
	10.217 17.335
	10.268 17.293
	10.338 17.342
	10.390 17.455
	10.411 17.521
	10.430 17.586
	10.465 17.685
	10.499 17.783
	10.517 17.847
	10.538 17.913
	10.590 18.027
	10.660 18.076
	10.711 18.033
	10.734 17.972
	10.757 17.879
	/
\linethickness=1pt
\setplotsymbol ({\makebox(0,0)[l]{\tencirc\symbol{'160}}})
%
% Fig INTERPOLATED PT SPLINE
%
\plot 11.443 18.076 	11.367 18.091
	11.301 18.102
	11.195 18.111
	11.114 18.102
	11.051 18.076
	10.978 17.997
	10.927 17.886
	10.889 17.773
	10.856 17.685
	10.822 17.586
	10.803 17.521
	10.782 17.455
	10.730 17.342
	10.660 17.293
	10.609 17.336
	10.586 17.396
	10.562 17.488
	/
\linethickness=1pt
\setplotsymbol ({\makebox(0,0)[l]{\tencirc\symbol{'160}}})
%
% Fig INTERPOLATED PT SPLINE
%
\plot 11.443 17.293 	11.367 17.278
	11.301 17.266
	11.195 17.258
	11.114 17.266
	11.051 17.293
	10.998 17.355
	10.976 17.410
	10.954 17.488
	/
%
% Fig CIRCULAR ARC object
%
\linethickness=1pt
\setplotsymbol ({\makebox(0,0)[l]{\tencirc\symbol{'160}}})
\circulararc 180.414 degrees from 11.932 22.430 center at 11.933 22.723
%
% Fig CIRCULAR ARC object
%
\linethickness=1pt
\setplotsymbol ({\makebox(0,0)[l]{\tencirc\symbol{'160}}})
\circulararc 180.000 degrees from 11.932 23.603 center at 11.932 23.897
%
% Fig POLYLINE object
%
\linethickness=1pt
\setplotsymbol ({\makebox(0,0)[l]{\tencirc\symbol{'160}}})
\putrule from  8.606 22.430 to  9.193 22.430
%
% Fig POLYLINE object
%
\linethickness=1pt
\setplotsymbol ({\makebox(0,0)[l]{\tencirc\symbol{'160}}})
\putrule from  9.193 22.333 to  9.193 23.132
\putrule from  9.193 23.114 to  9.993 23.114
\putrule from  9.976 23.114 to  9.976 22.315
\putrule from  9.976 22.333 to  9.175 22.333
%
% Fig POLYLINE object
%
\linethickness=1pt
\setplotsymbol ({\makebox(0,0)[l]{\tencirc\symbol{'160}}})
\putrule from  8.606 23.017 to  9.193 23.017
%
% Fig POLYLINE object
%
\linethickness=1pt
\setplotsymbol ({\makebox(0,0)[l]{\tencirc\symbol{'160}}})
\putrule from  8.606 23.603 to 10.562 23.603
%
% Fig POLYLINE object
%
\linethickness=1pt
\setplotsymbol ({\makebox(0,0)[l]{\tencirc\symbol{'160}}})
\putrule from  8.606 24.191 to 11.932 24.191
%
% Fig POLYLINE object
%
\linethickness=1pt
\setplotsymbol ({\makebox(0,0)[l]{\tencirc\symbol{'160}}})
\putrule from 11.345 23.603 to 11.932 23.603
%
% Fig POLYLINE object
%
\linethickness=1pt
\setplotsymbol ({\makebox(0,0)[l]{\tencirc\symbol{'160}}})
\putrule from  9.976 23.017 to 10.562 23.017
%
% Fig POLYLINE object
%
\linethickness=1pt
\setplotsymbol ({\makebox(0,0)[l]{\tencirc\symbol{'160}}})
\putrule from 11.345 23.017 to 11.932 23.017
%
% Fig POLYLINE object
%
\linethickness=1pt
\setplotsymbol ({\makebox(0,0)[l]{\tencirc\symbol{'160}}})
\putrule from  9.976 22.430 to 11.932 22.430
%
% Fig POLYLINE object
%
\linethickness=1pt
\setplotsymbol ({\makebox(0,0)[l]{\tencirc\symbol{'160}}})
\putrule from 10.562 22.919 to 10.562 23.720
\putrule from 10.562 23.702 to 11.363 23.702
\putrule from 11.345 23.702 to 11.345 22.902
\putrule from 11.345 22.919 to 10.545 22.919
%
% Fig CIRCULAR ARC object
%
\linethickness=1pt
\setplotsymbol ({\makebox(0,0)[l]{\tencirc\symbol{'160}}})
\circulararc 179.586 degrees from  3.090 23.017 center at  3.091 22.723
%
% Fig CIRCULAR ARC object
%
\linethickness=1pt
\setplotsymbol ({\makebox(0,0)[l]{\tencirc\symbol{'160}}})
\circulararc 179.172 degrees from  3.090 24.191 center at  3.092 23.897
%
% Fig POLYLINE object
%
\linethickness=1pt
\setplotsymbol ({\makebox(0,0)[l]{\tencirc\symbol{'160}}})
\putrule from  3.679 22.919 to  3.679 23.720
\putrule from  3.679 23.702 to  4.477 23.702
\putrule from  4.460 23.702 to  4.460 22.902
\putrule from  4.460 22.919 to  3.661 22.919
%
% Fig POLYLINE object
%
\linethickness=1pt
\setplotsymbol ({\makebox(0,0)[l]{\tencirc\symbol{'160}}})
\putrule from  5.048 22.333 to  5.048 23.132
\putrule from  5.048 23.114 to  5.847 23.114
\putrule from  5.829 23.114 to  5.829 22.315
\putrule from  5.829 22.333 to  5.031 22.333
%
% Fig POLYLINE object
%
\linethickness=1pt
\setplotsymbol ({\makebox(0,0)[l]{\tencirc\symbol{'160}}})
\putrule from  6.418 22.919 to  6.418 23.720
\putrule from  6.418 23.702 to  7.219 23.702
\putrule from  7.201 23.702 to  7.201 22.902
\putrule from  7.201 22.919 to  6.400 22.919
%
% Fig POLYLINE object
%
\linethickness=1pt
\setplotsymbol ({\makebox(0,0)[l]{\tencirc\symbol{'160}}})
\putrule from  3.090 22.430 to  5.048 22.430
%
% Fig POLYLINE object
%
\linethickness=1pt
\setplotsymbol ({\makebox(0,0)[l]{\tencirc\symbol{'160}}})
\putrule from  3.090 23.017 to  3.679 23.017
%
% Fig POLYLINE object
%
\linethickness=1pt
\setplotsymbol ({\makebox(0,0)[l]{\tencirc\symbol{'160}}})
\putrule from  3.090 23.603 to  3.679 23.603
%
% Fig POLYLINE object
%
\linethickness=1pt
\setplotsymbol ({\makebox(0,0)[l]{\tencirc\symbol{'160}}})
\putrule from  3.090 24.191 to  7.787 24.191
%
% Fig POLYLINE object
%
\linethickness=1pt
\setplotsymbol ({\makebox(0,0)[l]{\tencirc\symbol{'160}}})
\putrule from  4.460 23.603 to  6.418 23.603
%
% Fig POLYLINE object
%
\linethickness=1pt
\setplotsymbol ({\makebox(0,0)[l]{\tencirc\symbol{'160}}})
\putrule from  4.460 23.017 to  5.048 23.017
%
% Fig POLYLINE object
%
\linethickness=1pt
\setplotsymbol ({\makebox(0,0)[l]{\tencirc\symbol{'160}}})
\putrule from  5.829 23.017 to  6.418 23.017
%
% Fig POLYLINE object
%
\linethickness=1pt
\setplotsymbol ({\makebox(0,0)[l]{\tencirc\symbol{'160}}})
\putrule from  5.829 22.430 to  7.787 22.430
%
% Fig POLYLINE object
%
\linethickness=1pt
\setplotsymbol ({\makebox(0,0)[l]{\tencirc\symbol{'160}}})
\putrule from  7.201 23.017 to  7.787 23.017
%
% Fig POLYLINE object
%
\linethickness=1pt
\setplotsymbol ({\makebox(0,0)[l]{\tencirc\symbol{'160}}})
\putrule from  7.201 23.603 to  7.787 23.603
%
% Fig CIRCULAR ARC object
%
\linethickness=1pt
\setplotsymbol ({\makebox(0,0)[l]{\tencirc\symbol{'160}}})
\circulararc 180.000 degrees from  2.932 19.844 center at  2.932 19.550
%
% Fig CIRCULAR ARC object
%
\linethickness=1pt
\setplotsymbol ({\makebox(0,0)[l]{\tencirc\symbol{'160}}})
\circulararc 180.414 degrees from  2.932 21.016 center at  2.931 20.723
%
% Fig POLYLINE object
%
\linethickness=1pt
\setplotsymbol ({\makebox(0,0)[l]{\tencirc\symbol{'160}}})
\putrule from  3.518 19.744 to  3.518 20.545
\putrule from  3.518 20.527 to  4.319 20.527
\putrule from  4.301 20.527 to  4.301 19.727
\putrule from  4.301 19.744 to  3.500 19.744
%
% Fig POLYLINE object
%
\linethickness=1pt
\setplotsymbol ({\makebox(0,0)[l]{\tencirc\symbol{'160}}})
\putrule from  4.887 19.158 to  4.887 19.959
\putrule from  4.887 19.941 to  5.688 19.941
\putrule from  5.671 19.941 to  5.671 19.140
\putrule from  5.671 19.158 to  4.870 19.158
%
% Fig POLYLINE object
%
\linethickness=1pt
\setplotsymbol ({\makebox(0,0)[l]{\tencirc\symbol{'160}}})
\putrule from  6.257 19.744 to  6.257 20.545
\putrule from  6.257 20.527 to  7.058 20.527
\putrule from  7.040 20.527 to  7.040 19.727
\putrule from  7.040 19.744 to  6.239 19.744
%
% Fig POLYLINE object
%
\linethickness=1pt
\setplotsymbol ({\makebox(0,0)[l]{\tencirc\symbol{'160}}})
\putrule from  2.932 19.255 to  4.887 19.255
%
% Fig POLYLINE object
%
\linethickness=1pt
\setplotsymbol ({\makebox(0,0)[l]{\tencirc\symbol{'160}}})
\putrule from  2.932 19.844 to  3.518 19.844
%
% Fig POLYLINE object
%
\linethickness=1pt
\setplotsymbol ({\makebox(0,0)[l]{\tencirc\symbol{'160}}})
\putrule from  2.932 20.430 to  3.518 20.430
%
% Fig POLYLINE object
%
\linethickness=1pt
\setplotsymbol ({\makebox(0,0)[l]{\tencirc\symbol{'160}}})
\putrule from  2.932 21.016 to  7.626 21.016
%
% Fig POLYLINE object
%
\linethickness=1pt
\setplotsymbol ({\makebox(0,0)[l]{\tencirc\symbol{'160}}})
\putrule from  4.301 20.430 to  6.257 20.430
%
% Fig POLYLINE object
%
\linethickness=1pt
\setplotsymbol ({\makebox(0,0)[l]{\tencirc\symbol{'160}}})
\putrule from  4.301 19.844 to  4.887 19.844
%
% Fig POLYLINE object
%
\linethickness=1pt
\setplotsymbol ({\makebox(0,0)[l]{\tencirc\symbol{'160}}})
\putrule from  5.671 19.844 to  6.257 19.844
%
% Fig POLYLINE object
%
\linethickness=1pt
\setplotsymbol ({\makebox(0,0)[l]{\tencirc\symbol{'160}}})
\putrule from  5.671 19.255 to  7.626 19.255
%
% Fig POLYLINE object
%
\linethickness=1pt
\setplotsymbol ({\makebox(0,0)[l]{\tencirc\symbol{'160}}})
\putrule from  7.040 19.844 to  7.626 19.844
%
% Fig POLYLINE object
%
\linethickness=1pt
\setplotsymbol ({\makebox(0,0)[l]{\tencirc\symbol{'160}}})
\putrule from  7.040 20.430 to  7.626 20.430
%
% Fig CIRCULAR ARC object
%
\linethickness=1pt
\setplotsymbol ({\makebox(0,0)[l]{\tencirc\symbol{'160}}})
\circulararc 179.586 degrees from 11.737 19.850 center at 11.736 20.143
%
% Fig CIRCULAR ARC object
%
\linethickness=1pt
\setplotsymbol ({\makebox(0,0)[l]{\tencirc\symbol{'160}}})
\circulararc 192.620 degrees from 11.737 19.262 center at 11.834 20.142
%
% Fig POLYLINE object
%
\linethickness=1pt
\setplotsymbol ({\makebox(0,0)[l]{\tencirc\symbol{'160}}})
\putrule from  8.410 20.436 to  8.996 20.436
%
% Fig POLYLINE object
%
\linethickness=1pt
\setplotsymbol ({\makebox(0,0)[l]{\tencirc\symbol{'160}}})
\putrule from  8.410 19.850 to  8.996 19.850
%
% Fig POLYLINE object
%
\linethickness=1pt
\setplotsymbol ({\makebox(0,0)[l]{\tencirc\symbol{'160}}})
\putrule from  8.996 20.534 to  8.996 19.733
\putrule from  8.996 19.751 to  9.797 19.751
\putrule from  9.779 19.751 to  9.779 20.551
\putrule from  9.779 20.534 to  8.978 20.534
%
% Fig POLYLINE object
%
\linethickness=1pt
\setplotsymbol ({\makebox(0,0)[l]{\tencirc\symbol{'160}}})
\putrule from 10.367 19.947 to 10.367 19.147
\putrule from 10.367 19.164 to 11.166 19.164
\putrule from 11.148 19.164 to 11.148 19.965
\putrule from 11.148 19.947 to 10.350 19.947
%
% Fig POLYLINE object
%
\linethickness=1pt
\setplotsymbol ({\makebox(0,0)[l]{\tencirc\symbol{'160}}})
\putrule from  9.779 19.850 to 10.367 19.850
%
% Fig POLYLINE object
%
\linethickness=1pt
\setplotsymbol ({\makebox(0,0)[l]{\tencirc\symbol{'160}}})
\putrule from  8.410 19.262 to 10.367 19.262
%
% Fig POLYLINE object
%
\linethickness=1pt
\setplotsymbol ({\makebox(0,0)[l]{\tencirc\symbol{'160}}})
\putrule from  8.410 21.023 to 11.737 21.023
%
% Fig POLYLINE object
%
\linethickness=1pt
\setplotsymbol ({\makebox(0,0)[l]{\tencirc\symbol{'160}}})
\putrule from  9.779 20.436 to 11.737 20.436
%
% Fig POLYLINE object
%
\linethickness=1pt
\setplotsymbol ({\makebox(0,0)[l]{\tencirc\symbol{'160}}})
\putrule from 11.148 19.850 to 11.737 19.850
%
% Fig POLYLINE object
%
\linethickness=1pt
\setplotsymbol ({\makebox(0,0)[l]{\tencirc\symbol{'160}}})
\putrule from 11.148 19.262 to 11.737 19.262
%
% Fig TEXT object
%
\put{$n$ even} [lB] at  7.117 18.578
% used to be 7.137
%
% Fig TEXT object
%
\put{$a_1$} [lB] at  3.715 20.045
%
% Fig TEXT object
%
\put{$a_2$} [lB] at  5.084 19.459
%
% Fig TEXT object
%
\put{$a_3$} [lB] at  6.454 20.045
%
% Fig TEXT object
%
\put{$a_{n-1}$} [lB] at  9.065 20.045
% used to be 9.095
%
% Fig TEXT object
%
\put{$a_n$} [lB] at 10.562 19.459
%
% Fig TEXT object
%
\put{$\ldots$} [lB] at  7.806 19.850
%
% Fig TEXT object
%
\put{$\ldots$} [lB] at  7.806 21.023
%
% Fig TEXT object
%
\put{$\ldots$} [lB] at  7.806 19.262
%
% Fig TEXT object
%
\put{$n$ odd} [lB] at  7.137 21.745
%
% Fig TEXT object
%
\put{$a_{n-1}$} [lB] at  9.280 22.625
% used to be 9.290
%
% Fig TEXT object
%
\put{$a_n$} [lB] at 10.757 23.213
%
% Fig TEXT object
%
\put{$a_3$} [lB] at  6.648 23.213
%
% Fig TEXT object
%
\put{$a_1$} [lB] at  3.909 23.213
%
% Fig TEXT object
%
\put{$a_2$} [lB] at  5.279 22.625
%
% Fig TEXT object
%
\put{$\ldots$} [lB] at  8.000 24.191
%
% Fig TEXT object
%
\put{$\ldots$} [lB] at  8.000 23.603
%
% Fig TEXT object
%
\put{$\ldots$} [lB] at  8.000 22.430
%
% Fig TEXT object
%
\put{$\ldots$} [lB] at  8.000 23.017
% used to be 8.079
%
% Fig TEXT object
%
\put{$\ldots$} [lB] at  7.826 20.436
% used to be 7.885
%
% Fig TEXT object
%
\put{positive twists} [lB] at  3.929 16.607
% used to be 3.909
%
% Fig TEXT object
%
\put{negative twists} [lB] at  9.193 16.607
\linethickness=0pt
\putrectangle corners at  2.601 24.238 and 12.753 16.546
\endpicture
}
\caption{The rational-form diagram $T(a_1,\dots,a_n)$.  
Each box contains the specified number
of half-twists; positive and negative twists are shown.
}
%\vspace{0.5in}
\label{bridge}
\end{figure}

\begin{figure}[ht!]
\centerline{
\includegraphics[width=2.6in,angle=270]{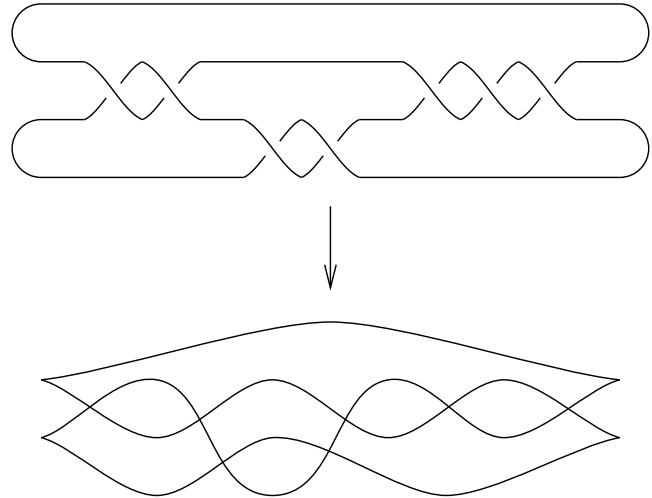}
}
% used to be width=2.8in
\vspace{0.2in}
\caption{The correspondence between a diagram in Legendrian rational form 
(in this case, $T(2,2,3)$, or $5_2$) 
and the front of a Legendrian link of the same ambient type.
}
\label{cusps}
\end{figure}

To each $T(a_1,\dots,a_n)$, we may associate a rational number (or $\infty$),
the continued fraction
\[
[a_1,\dots,a_n] = a_1 + \frac{1}{-a_2} \!\!
\begin{array}{l}
\\ + \end{array}
\!\!
\frac{1}{a_3} \!\!
\begin{array}{l}
\\ + \end{array}
\!\!
\frac{1}{-a_4} \!\!
\begin{array}{l}
\\ + \, \cdots \, + \end{array}
\!\!
\frac{1}{(-1)^{n-1} a_n}\,.
\]
Note that our convention is the opposite of the convention in
\cite{Lic}, and differs by alternating signs from the standard 
convention from, e.g., \cite{BZ}.
The classification of $2$-bridge links further states that
if $1/[a_1,\ldots,a_n] - 1/[b_1,\ldots,b_n] \in \Z$, then 
$T(a_1,\dots,a_n)$ and $T(b_1,\ldots,b_n)$ are ambient isotopic.
(The criterion stating precisely when two such links are isotopic is
only slightly more complicated, but will not concern us here.)

Now define $T(a_1,\dots,a_n)$ to be in {\it Legendrian rational form} if
$a_i \geq 2$ for all $i$.  
Although $T(a_1,\ldots,a_n)$ corresponds to a Legendrian link
whenever $a_i \geq 1$ for all $i$, it is crucial to the proof of
Lemma~\ref{lickorish} below that $a_i \geq 2$ for $2 \leq i \leq n-1$.  Indeed,
if one of these $a_i$ is $1$, then it is straightforward to see, by
drawing the front, that the resulting Legendrian link does not maximize
Thurston-Bennequin number.

Any link diagram in Legendrian rational 
form is easily
converted into the front (i.e., projection to the $xz$ plane) of a Legendrian
link by replacing the four vertical tangencies by cusps; see 
Figure~\ref{cusps}.  Since the crossings in a front are resolved locally
so that the strand with more negative slope always lies over the strand
with more positive slope, a link diagram in Legendrian rational form
is ambient isotopic to the corresponding front.  (This observation
explains our choice of convention for positive versus negative twists.)

\begin{lemma}
Any $2$-bridge link can be expressed as a diagram in Legendrian 
rational form.
\end{lemma}

\begin{proof}
Let $K$ be a $2$-bridge link; let $T(a_1,\dots,a_n)$ be a rational-form
diagram for $K$, and write $[a_1,\dots,a_n] = p/q$ for $p,q\in\Z$.
The classification of $2$-bridge links implies that $K$ is isotopic to
any rational-form diagram associated to the fraction
$r = p/(q-\lfloor \frac{q}{p} \rfloor p) > 1$.  (If $q/p$ is an integer,
then it is easy to see that $K$ is the trivial knot, which is not
$2$-bridge.)

Define a sequence $x_1,x_2,\ldots$ of rational numbers by $x_1 = r$,
$x_{i+1} = 1/(\lceil x_i \rceil - x_i)$.  This sequence terminates at,
say, $x_m$, where $x_m$ is an integer.  Write $b_i = \lceil x_i \rceil$.
It is easy to see that $b_i \geq 2$ for all $i$, and that
$r = [b_1,\dots,b_m]$.  Then $K$ is isotopic to
$T(b_1,\dots,b_m)$, which is in Legendrian rational form.
\end{proof}

Consider a link diagram $T=T(a_1,\dots,a_n)$ in Legendrian rational form,
and let $K$ be the (Legendrian link given by the) front obtained 
from $T$.  We claim that the Thurston-Bennequin number of $K$ 
agrees precisely with the Kauffman bound.  Recall that the Kauffman polynomial
$F_T(a,x)$ of $T$ is $a^{w(T)}$ times the {\it unoriented} Kauffman polynomial
(or L-polynomial) $L_T(a,x)$, where $w(T)$ is the writhe of the diagram $T$.
(Here we use Kauffman's original notation \cite{Kau}, except with $a$ 
replaced by $1/a$.)

We will need a matrix formula for $L_T(a,x)$ due to \cite{Lic}.
Write 
\begin{align*}
M = \left( \begin{smallmatrix}
x & -1 & x \\
1 & 0 & 0 \\
0 & 0 & 1/a
\end{smallmatrix} \right),~~~
&&
S = \left( \begin{smallmatrix}
0 & 1 & 0 \\
0 & 0 & 1 \\
1/a & 0 & 0
\end{smallmatrix} \right),~~~
&&
v = \left( \begin{smallmatrix}
1 \\ 0 \\ 0
\end{smallmatrix} \right),~~~
&&
w = \left( \begin{smallmatrix}
a \\ a^2 \\ \textstyle{\frac{a^2+1}{x}-a}
\end{smallmatrix} \right);
\end{align*}
\noindent then
\[
L_{T(a_1,\dots,a_n)}(a,x) = \textstyle{\frac{1}{a}} 
v^t M^{-a_1-1} S M^{-a_2-1} S \cdots
M^{-a_n-1} S w,
\]
where $t$ denotes transpose.

\begin{lemma}
If $a_1,a_n \geq 1$ and $a_i \geq 2$ for $2 \leq i \leq n-1$, then we have
$\mindeg L_{T(a_1,\dots,a_n)}(a,x) = -1$.
\label{lickorish}
\end{lemma}

\begin{proof}
None of $M^{-1}$, $M^{-1}S$, and $w$ contains negative powers of $a$;
the lemma will be proved if we can show that $f(x) \neq 0$, where
\[
f(x) = \left(v^t M^{-a_1} (M^{-1}S) M^{-a_2} (M^{-1}S) \cdots
M^{-a_n} (M^{-1}S) w\right) |_{a=0}.
\]
Define the auxiliary matrices 
\begin{align*}
A = M^{-1}|_{a=0} = \left( \begin{smallmatrix}
0 & 1 & 0 \\
-1 & x & 0 \\
0 & 0 & 0
\end{smallmatrix} \right),
&&
B = (M^{-2} S M^{-1})|_{a=0} = \left( \begin{smallmatrix}
1 & 0 & 0 \\
x & 0 & 0 \\
0 & 0 & 0
\end{smallmatrix} \right),
&&
u = \left( \begin{smallmatrix}
0 \\ 1 \\ 0 \end{smallmatrix} \right).
\end{align*}
%Define the auxiliary matrices $A = M^{-1}|_{a=0} = \left( \begin{smallmatrix}
%0 & 1 & 0 \\
%-1 & x & 0 \\
%0 & 0 & 0
%\end{smallmatrix} \right)$,
%$B = (M^{-2} S M^{-1})|_{a=0} = \left( \begin{smallmatrix}
%1 & 0 & 0 \\
%x & 0 & 0 \\
%0 & 0 & 0
%\end{smallmatrix} \right)$, and
%$u = (0,1,0)$.  
Then $(ASw)|_{a=0} = \frac{1}{x} Au$ and
$B = Au v^t$, and so
\begin{eqnarray*}
f(x) &=& v^t A^{a_1-1} B A^{a_2-2} B A^{a_3-2} B \cdots A^{a_{n-1}-2} B
A^{a_n-1} (ASw)|_{a=0} \\
&=& \textstyle{\frac{1}{x}}
(v^t A^{a_1} u)(v^t A^{a_2-1} u)(v^t A^{a_3-1} u) \cdots
(v^t A^{a_{n-1}-1} u) (v^t A^{a_n} u).
\end{eqnarray*}
But if we define a sequence of functions $f_k(x) = v^t A^k u$, 
then an easy induction
yields the recursion $f_{k+2}(x) = x f_{k+1}(x) - f_k(x)$ with 
$f_1(x) = 1$ and $f_2(x) = x$.  In particular,
for all $k \geq 1$, $f_k(x)$ has degree $k-1$ and is thus nonzero.
From the given conditions on $a_i$, it follows that $f(x) \neq 0$, as desired.
\end{proof}

%\vspace{16pt}

%\renewcommand{\proofname}{Proof of Theorem \ref{main}}
\begin{proof}[Proof of Theorem \ref{main}]
Let $T$ be a Legendrian rational form for a $2$-bridge link $K$.
The crossings of $T$ are counted, with the same signs, by both the
writhe of $T$ and the Thurston-Bennequin number of the Legendrian
link $K'$ obtained from $T$; $tb\:\!(K')$, however, also
subtracts half the number of cusps.  Hence
\begin{eqnarray*}
tb\:\!(K') &=& w(T) - 2 \\
&=& (\mindeg F_T(a,x) - \mindeg L_T(a,x)) - 2 \\
&=& \mindeg F_T(a,x) - 1
\end{eqnarray*}
by Lemma~\ref{lickorish}.
Since $K'$ is ambient isotopic to $K$, we conclude that
$\mtb{K}$ is at least $\mindeg F_T(a,x) - 1$; by the Kauffman bound, 
equality must hold.
\end{proof}

%\renewcommand{\proofname}{Proof}

%\newpage

\section*{Appendix: Maximal Thurston-Bennequin number for small knots}

The following table gives the maximal Thurston-Bennequin invariant
for all prime knots with nine or fewer crossings.  We distinguish between
mirrors by using the diagrams in \cite{Rol}: the knots $K$ are the
ones drawn in \cite{Rol}, with mirrors $\tilde{K}$.  A dagger next to
a knot indicates that it is not two-bridge; a double dagger indicates that
the knot is amphicheiral (identical to its unoriented mirror).
For the interested reader, two-bridge descriptions of the
two-bridge knots in the table can be deduced from the tables in
\cite{BZ}.

The boldfaced numbers indicate the knots for which the Kauffman
bound is not sharp (for $8_{19}$), or probably not sharp (for $9_{42}$).  
As mentioned in the Introduction, it is believed
that $\tb=-5$ for the mirror $9_{42}$ knot; the best known bound, however,
is the Kauffman bound $\tb \leq -3$.

\[
\begin{array}{|l|r|r||l|r|r||l|r|r|}
\hline
\rule{0pt}{14pt} K & \mtb{K} & \mtb{\tilde{K}} &
K & \mtb{K} & \mtb{\tilde{K}} &
K & \mtb{K} & \mtb{\tilde{K}} \\
\hline\hline
0_1 	& -1 	& \ddagger 	&
{8_{15}}^\dagger	&-13	&3	&
{9_{22}}^\dagger	&-3	&-8     \\
3_1	&-6	&1	&
{8_{16}}^\dagger	&-8	&-2	&
9_{23}	&-14	&3	\\
4_1	&-3	&\ddagger	&
{8_{17}}^\dagger	&-5	&\ddagger	&
{9_{24}}^\dagger	&-6	&-5	\\
5_1	&-10	&3	&
{8_{18}}^\dagger	&-5	&\ddagger	&
{9_{25}}^\dagger	&-10	&-1	\\
5_2	&-8	&1	&
{8_{19}}^\dagger	&5	&\mathbf{-12}	&
9_{26}	&-2	&-9	\\
6_1	&-5	&-3	&
{8_{20}}^\dagger	&-6	&-2	&
9_{27}	&-6	&-5	\\
6_2	&-7	&-1	&
{8_{21}}^\dagger	&-9	&1	&
{9_{28}}^\dagger	&-9	&-2	\\
6_3	&-4	&\ddagger	&
9_1	&-18	&7	&
{9_{29}}^\dagger 	& -8 	& -3 	\\
7_1	&-14	&5	&
9_2	&-12	&1	&
{9_{30}}^\dagger	&-6	&-5	\\
7_2	&-10	&1	&
9_3	&5	&-16	&
9_{31}	&-9	&-2	\\
7_3	&3	&-12	&
9_4	&-14	&3	&
{9_{32}}^\dagger	&-2	&-9	\\
7_4	&1	&-10	&
9_5	&1	&-12	&
{9_{33}}^\dagger	&-6	&-5	\\
7_5	&-12	&3	&
9_6	&-16	&5	&
{9_{34}}^\dagger	&-6	&-5	\\
7_6	&-8	&-1	&
9_7	&-14	&3	&
{9_{35}}^\dagger	&-12	&1	\\
7_7	&-4	&-5	&
9_8 	& -8 	& -3 	&
{9_{36}}^\dagger	&1	&-12	\\
8_1	&-7	&-3	&
9_9	&-16	&5	&
{9_{37}}^\dagger	&-6	&-5	\\
8_2	&-11	&1	&
9_{10}	&3	&-14	&
{9_{38}}^\dagger	&-14	&3	\\
8_3	&-5	&\ddagger	&
9_{11}	&1	&-12	&
{9_{39}}^\dagger	&-1	&-10	\\
8_4	&-7	&-3	&
9_{12}	&-10	&-1	&
{9_{40}}^\dagger	&-9	&-2	\\
{8_5}^\dagger	&1	&-11	&
9_{13}	&3	&-14	&
{9_{41}}^\dagger	&-7	&-4	\\
8_6	&-9	&-1	&
9_{14}	&-4	&-7	&
{9_{42}}^\dagger	&-3	&\mathbf{-5} (?)	\\
8_7	&-2	&-8	&
9_{15}	&-10	&-1	&
{9_{43}}^\dagger	&1	&-10	\\
8_8 	& -4 	& -6 	&
{9_{16}}^\dagger	&5	&-16	&
{9_{44}}^\dagger	&-6	&-3	\\
8_9	&-5	&\ddagger	&
9_{17}	&-8	&-3	&
{9_{45}}^\dagger	&-10	&1	\\
{8_{10}}^\dagger	&-2	&-8	&
9_{18}	&-14	&3	&
{9_{46}}^\dagger	&-7	&-1	\\
8_{11}	&-9	&-1	&
9_{19}	&-6	&-5	&
{9_{47}}^\dagger	&-2	&-7	\\
8_{12}	&-5	&\ddagger	&
9_{20}	&-12	&1	&
{9_{48}}^\dagger	&-1	&-8	\\
8_{13}	&-4	&-6	&
9_{21}	&-1	&-10	&
{9_{49}}^\dagger	&3	&-12	\\
8_{14}	&-9	&-1	&
&&&
&& \\
\hline
\end{array}
\]

\Addresses\recd


\newpage

\begin{thebibliography}

\bibitem{Ben}
{\bf D Bennequin}, {\it Entrelacements et \'equations de Pfaff}, 
Ast\'erisque
107--108 (1983) 87--161

\bibitem{BZ}
{\bf G Burde}, {\bf H Zieschang}, {\it Knots}, Walter de Gruyter
(1985)

\bibitem{Eps}
{\bf J Epstein}, {\it On the invariants and isotopies of 
Legendrian and transverse
knots}, PhD thesis, U\,C Davis, 1997

\bibitem{EH}
{\bf J Etnyre}, {\bf K Honda}, {\it Knots and contact geometry},
preprint, 2000, {\tt arXiv: math.GT/0006112}

\bibitem{Fer}
{\bf E Ferrand}, {\it On Legendrian knots and polynomial invariants},
preprint, 2000, {\tt arXiv:math.GT/0002250}

\bibitem{FT}
{\bf D Fuchs}, {\bf S Tabachnikov}, {\it Invariants of Legendrian 
and transverse
knots in the standard contact space}, Topology 36 (1997) 1025--1054

\bibitem{Kau}
{\bf L Kauffman}, {\it On Knots}, Princeton University Press (1987)

\bibitem{Lic}
{\bf W\,B\,R Lickorish}, {\it Linear skein theory and link polynomials},
Topology Appl.\ 27 (1987) 265--274

%\bibitem{Mur}
%{\bf K Murasugi}, {\it Knot Theory and Its Applications}, Birkh\"auser
%(1996)

\bibitem{Ng}
{\bf L Ng}, {\it Invariants of Legendrian links}, PhD thesis, MIT, 2001,
available at
{\tt http://www-math.mit.edu/\~{}lenny/}

\bibitem{Rol}
{\bf D Rolfsen}, {\it Knots and Links}, Publish or Perish (1990)

\bibitem{Rud}
{\bf L Rudolph}, {\it A congruence between link polynomials},
Math.\ Proc.\ Camb.\ Phil.\ Soc.\ 107 (1990) 319--327

\bibitem{Sch}
{\bf H Schubert}, {\it Knoten mit zwei Br\"ucken}, Math.\ Zeit. 65 (1956)
133--170

%\bibitem[Tab]{Tab}
%S.\ Tabachnikov, Bennequin number of Legendrian links, {\it Math.\
%Res.\ Lett.} {\bf 4} (1997), 143--156.

\bibitem{Tan}
{\bf T Tanaka}, {\it Maximal Bennequin numbers and Kauffman polynomials of 
positive links}, Proc.\ Amer.\ Math.\ Soc.\ 127
(1999) 3427--3432

\bibitem{Yu}
{\bf N Yufa}, {\it Thurston-Bennequin invariant
of Legendrian knots}, senior thesis, MIT, 2001

\end{thebibliography}
\end{document}